\documentclass[a4paper,12pt,reqno]{amsart}

\usepackage{amsmath,amssymb,amsthm,mathrsfs}
\usepackage{mathtools}
\usepackage[all,cmtip]{xy}
\usepackage{bm}
\usepackage{afterpage}
\usepackage[a4paper,margin=2cm]{geometry}
\usepackage{pdflscape}

\usepackage[shortlabels]{enumitem}

\newtheorem{theorem}{Theorem}[section]
\newtheorem{proposition}[theorem]{Proposition}

\newtheorem{corollary}[theorem]{Corollary}
\theoremstyle{definition} 
\newtheorem{definition}[theorem]{Definition} 
\newtheorem{remark}[theorem]{Remark}
\newtheorem{example}[theorem]{Example}

\makeatletter
\let\c@equation\c@theorem
\makeatother
\numberwithin{equation}{section}

\usepackage[hidelinks]{hyperref}

\newcommand{\lra}{\longrightarrow}
\newcommand{\Cat}{\operatorname{\mathbf{Cat}}}
\newcommand{\vcat}{\mathcal{V}\text{-}\mathbf{Cat}}
\newcommand{\vact}{\mathcal{V}\text{-}\mathbf{Act}}
\newcommand{\vskcat}{\mathcal{V}\text{-}\mathbf{SkCat}}
\newcommand{\vskact}{\mathcal{V}\text{-}\mathbf{SkAct}}
\newcommand{\vskproact}{\mathcal{V}\text{-}\mathbf{SkProAct}}

\newcommand{\skmoncat}{\operatorname{\mathbf{SkMonCat}}}
\newcommand{\skenrcat}{\operatorname{\mathbf{SkEnrCat}}}
\newcommand{\twocat}{\mathbf{2}\text{-}\mathbf{Cat}}
\newcommand{\bigtwocat}{\mathbf{2}\text{-}\mathbf{CAT}}
\newcommand{\Set}{\operatorname{\mathbf{Set}}}
\newcommand{\Hom}{\mathrm{Hom}}

\newcommand{\cd}[2][]{\vcenter{\hbox{\xymatrix#1{#2}}}}

\makeatletter
\def\matrixobject@{%
  \edef \next@{={\DirectionfromtheDirection@ }}%
  \expandafter \toks@ \next@ \plainxy@
  \let\xy@@ix@=\xyq@@toksix@
  \xyFN@ \OBJECT@}
\let\xy@entry@@norm=\entry@@norm
\def\entry@@norm@patched{%
  \let\object@=\matrixobject@
  \xy@entry@@norm }
\AtBeginDocument{\let\entry@@norm\entry@@norm@patched}
\makeatother

\newcommand{\hdash}{\rotatebox[origin=c]{90}{$\vdash$}}

\title{Skew-enriched categories}
\author{Alexander Campbell}
\address{Centre of Australian Category Theory, Department of Mathematics, Macquarie University, NSW 2109, Australia}
\email{alexander.campbell@mq.edu.au}

\subjclass[2010]{18D20, 18D10, 18D15}
\date{October 6, 2017}

\begin{document}
\maketitle

\begin{abstract}
This paper introduces a skew variant of the notion of enriched category, suitable for enrichment over a skew-monoidal category, the main novelty of which is that the elements of the enriched hom-objects need not be in bijection with the morphisms of the underlying category. This is the natural setting in which to introduce the notion of locally weak comonad, which is fundamental to the theory of enriched algebraic weak factorisation systems. The equivalence, for a monoidal closed category $\mathcal{V}$, between tensored $\mathcal{V}$-categories and hommed $\mathcal{V}$-actegories is extended to the skew setting and easily proved by recognising both skew $\mathcal{V}$-categories and skew $\mathcal{V}$-actegories as equivalent to special kinds of skew $\mathcal{V}$-proactegory. 
\end{abstract}

\tableofcontents

\section{Introduction}
Recent years have seen the introduction of ``skew'' variants of the categorical notions of monoidal category \cite{MR2964621}, closed category \cite{MR3010098}, and actegory \cite{MR3370862}. The definitions of these skew structures differ from those of their classical counterparts by not demanding the invertibility  of certain structural morphisms, for which directions must therefore be specified. 
For example, the associativity and unit constraints of a (left) skew-monoidal category are natural transformations $a \colon (X \otimes Y) \otimes Z \lra X \otimes (Y\otimes Z)$, $l \colon I \otimes X \lra X$, and $r \colon X \lra X \otimes I$, none of which need be invertible; for a skew-closed category $\mathcal{V}$, the natural transformation $i \colon [I,X] \lra X$ and the composite functions
\begin{equation} \label{wkhom}
\cd[@C=3em]{
\mathcal{V}(X,Y) \ar[r]^-{[X,-]} & \mathcal{V}([X,X],[X,Y]) \ar[r]^-{\mathcal{V}(j,1)} & \mathcal{V}(I,[X,Y])
}
\end{equation}
are not required to be invertible. 

The extra freedom thus granted is well illustrated in the examples from \cite{Bourke2015} of skew monoidal closed structures on Quillen model categories  (which descend to genuine monoidal closed structures on the homotopy categories), where the elements of the domain and codomain of the function (\ref{wkhom}) can be interpreted as ``strict morphisms'' and ``weak morphisms'' respectively; for instance, there is a skew monoidal closed category of permutative categories, whose strict and weak morphisms are the symmetric strict monoidal functors and symmetric strong monoidal functors respectively.   A perhaps more familiar example is the skew monoidal closed structure on the category $\textbf{Gp}$ of groups and group homomorphisms, whose internal hom $[G,H]$ is the set of arbitrary functions from $G$ to $H$ equipped with the pointwise group structure. Furthermore, whereas the category $\Set$ of sets and functions admits only one biclosed monoidal structure, each monoid $M$ determines a biclosed skew-monoidal structure on $\Set$ with tensor product $X \otimes Y = M \times X \times Y$.

This paper introduces a skew variant of the notion of enriched category,  the main novelty of which is that the elements of the enriched hom-objects need not be in bijection with the morphisms of the underlying category. The formulation of this notion requires that a skew-enriched category be conceived as a category with structure, and not as an independent structure from which an underlying category is derived, as in the common definition of enriched category (e.g.\ \cite{MR2177301}). Thus, for a skew-monoidal category $\mathcal{V}$, a skew $\mathcal{V}$-category has hom-sets $\mathcal{A}(A,B)$ and enriched hom-objects $\underline{\mathcal{A}}(A,B)$ belonging to $\mathcal{V}$, and the canonical functions
\begin{equation} \label{strwk0}
\cd[@C=3em]{
\mathcal{A}(A,B) \ar[r]^-{\underline{\mathcal{A}}(A,-)} & \mathcal{V}(\underline{\mathcal{A}}(A,A),\underline{\mathcal{A}}(A,B)) \ar[r]^-{\mathcal{V}(j,1)} & \mathcal{V}(I,\underline{\mathcal{A}}(A,B))
}
\end{equation}
need not be invertible.     In Section \ref{skewenrsection} we define the $2$-category of skew $\mathcal{V}$-categories, and show that its full sub-$2$-category on those skew $\mathcal{V}$-categories for which the functions (\ref{strwk0}) are invertible is equivalent to the $2$-category of $\mathcal{V}$-categories in the usual sense (defined for a skew-monoidal $\mathcal{V}$ in \cite{MR3010098}), at least when the left unit constraint $l$ of the base $\mathcal{V}$ is invertible.

Every skew monoidal closed category is skew-enriched over itself, and the interpretation of the elements of the domain and codomain of (\ref{strwk0}) as strict morphisms and weak morphisms respectively may be profitably carried over to this more general context. For example, the category of groups in a (locally small) category $\mathscr{C}$ with finite products is skew-enriched over the skew monoidal closed structure on $\textbf{Gp}$, with strict and weak morphisms given by the group homomorphisms and arbitrary morphisms in $\mathscr{C}$ respectively. For a monoid $M$, each $M$-set $A$ determines a skew-enrichment of any locally small category $\mathscr{E}$ over the skew-monoidal structure on $\Set$ determined by $M$, with enriched hom-sets given by the powers $\mathscr{E}(X,Y)^A$.

Our main motivation for introducing the notion of skew-enriched category comes from a problem in enriched homotopy theory. It is known that if $\mathcal{V}$ is a monoidal model category in which every object is cofibrant, then any cofibrantly generated $\mathcal{V}$-enriched model category has a $\mathcal{V}$-enriched cofibrant replacement comonad (see \cite{MR3221774}); conversely, if a monoidal model category $\mathcal{V}$ (with cofibrant unit object) has a $\mathcal{V}$-enriched cofibrant replacement comonad, then every object of $\mathcal{V}$ must be cofibrant \cite{MR3531994}. These results leave open the question of what extra structure, if not an enrichment in the ordinary sense, is naturally possessed by the cofibrant replacement comonad of an enriched model category when not every object of the base monoidal model category is cofibrant. In a future paper (in the meantime see \cite{ct2017talk}) we will answer this question by showing that the cofibrant replacement comonad of an enriched algebraic weak factorisation system is a \emph{locally weak comonad}, a notion which we define in Section \ref{locwksection}.
    
The notion of skew-enriched category arose in the study of locally weak comonads, in particular as the natural context in which to state and prove Proposition \ref{locwkthm}, which shows how a new skew-enriched category can be constructed from a locally weak comonad on an old one. In our motivating examples, this shows that, for a monoidal algebraic weak factorisation system $(L,R)$ on a monoidal category $\mathcal{V}$ and an $(L,R)$-enriched algebraic weak factorisation system $(H,M)$ on a $\mathcal{V}$-category $\mathcal{A}$, there is a skew-enrichment of $\mathcal{A}$ over the skew-monoidal structure on $\mathcal{V}$ induced by the cofibrant replacement monoidal comonad for $(L,R)$, for which the weak morphisms are the weak maps for $(H,M)$ (in the sense of \cite{MR3393454}).

In Section \ref{proactsection} we further justify the notion of skew-enriched category by showing, for a skew-monoidal category $\mathcal{V}$, that to give a skew $\mathcal{V}$-enrichment of a category $\mathcal{A}$ such that each functor $\underline{\mathcal{A}}(A,-) \colon \mathcal{A}\lra \mathcal{V}$ has a left adjoint is equivalently to give a skew $\mathcal{V}$-action on $\mathcal{A}$ such that each functor $-\ast A \colon \mathcal{V} \lra \mathcal{A}$ has a right adjoint, thereby extending to the skew setting the equivalence, for a monoidal closed category $\mathcal{V}$, between tensored $\mathcal{V}$-categories and hommed $\mathcal{V}$-actegories. This follows immediately from the observation that both skew $\mathcal{V}$-categories and skew $\mathcal{V}$-actegories are equivalent to special kinds of skew $\mathcal{V}$-proactegory.

\subsection*{Acknowledgements}
The support of Australian Research Council Future Fellowship \linebreak FT160100393 is gratefully acknowledged. The author also thanks Richard Garner for his helpful comments on early presentations of this paper.

\section{Skew-enriched categories} \label{skewenrsection}
Let $(\mathcal{V},\otimes,I,a,l,r)$ be a (left) skew-monoidal category.   In this section we define the 2-category of skew $\mathcal{V}$-categories, and compare it with the 2-category of $\mathcal{V}$-categories in the usual sense.

\begin{definition}
A \emph{skew $\mathcal{V}$-category} $\underline{\mathcal{A}}$ consists of the following data:
\begin{enumerate}[(i)]
\item a category $\mathcal{A}$,
\item a functor $\underline{\mathcal{A}}(-,-) \colon \mathcal{A}^\mathrm{op} \times \mathcal{A} \lra \mathcal{V}$,
\item a natural transformation $M \colon \underline{\mathcal{A}}(B,C) \otimes \underline{\mathcal{A}}(A,B) \lra \underline{\mathcal{A}}(A,C)$,
\item a natural transformation $j  \colon I \lra \underline{\mathcal{A}}(A,A)$,
\end{enumerate}
subject to the following three axioms.
\begin{equation*}
\cd{
(\underline{\mathcal{A}}(C,D)\otimes\underline{\mathcal{A}}(B,C))\otimes\underline{\mathcal{A}}(A,B) \ar[r]^-{M\otimes 1} \ar[d]_-a  & \underline{\mathcal{A}}(B,D)\otimes\underline{\mathcal{A}}(A,B) \ar[dd]^-M \\
\underline{\mathcal{A}}(C,D) \otimes (\underline{\mathcal{A}}(B,C)\otimes\underline{\mathcal{A}}(A,B)) \ar[d]_-{1 \otimes M} & {} \\
\underline{\mathcal{A}}(C,D) \otimes \underline{\mathcal{A}}(A,C) \ar[r]_-M & \underline{\mathcal{A}}(A,D)
}
\end{equation*}
\begin{equation*}
\cd[@C=0em]{
I \otimes \underline{\mathcal{A}}(A,B) \ar[dr]_-l \ar[rr]^-{j\otimes 1} && \underline{\mathcal{A}}(B,B) \otimes \underline{\mathcal{A}}(A,B) \ar[dl]^-M \\
& \underline{\mathcal{A}}(A,B)
}
\qquad
\cd{
\underline{\mathcal{A}}(A,B) \ar[r]^-1 \ar[d]_-r & \underline{\mathcal{A}}(A,B) \\
\underline{\mathcal{A}}(A,B) \otimes I \ar[r]_-{1 \otimes j} & \underline{\mathcal{A}}(A,B) \otimes \underline{\mathcal{A}}(A,A) \ar[u]_-M
}
\end{equation*}

A skew $\mathcal{V}$-category $\underline{\mathcal{A}}$ is said to be \emph{normal} when the composite function
\begin{equation} \label{strwk}
\cd[@C=3em]{
\mathcal{A}(A,B) \ar[r]^-{\underline{\mathcal{A}}(A,-)} & \mathcal{V}(\underline{\mathcal{A}}(A,A),\underline{\mathcal{A}}(A,B)) \ar[r]^-{\mathcal{V}(j_A,1)} & \mathcal{V}(I,\underline{\mathcal{A}}(A,B))
}
\end{equation}
is a  bijection for each $A,B \in \mathcal{A}$.
\end{definition}

Recall that a skew-monoidal category is said to be \emph{left normal} when its left unit constraints $l \colon I \otimes X \lra X$ are invertible \cite{MR2972969}.

\begin{remark}
This definition departs from the usual definition of enriched category \cite{MR0225841,MR2177301}  by defining a (skew-)enriched category to be a structure on a given category $\mathcal{A}$, rather than an independent structure from which an underlying category is derived. Accordingly, in our definition we require the enriched hom-objects to be the values of a functor and the families of morphisms $M$ and $j$ to be natural with respect to the morphisms of $\mathcal{A}$. 

This departure frees us to capture the phenomenon of skew-enrichment, in which the elements of the enriched hom-objects need not be in bijection with the morphisms of the category $\mathcal{A}$. Note that the usual definition of $\mathcal{V}$-category carries over to a skew-monoidal base $\mathcal{V}$ (see \cite{MR3010098}); over a left normal base, $\mathcal{V}$-categories in the usual sense are equivalent to normal skew $\mathcal{V}$-categories (Corollary \ref{2adj}), whereas over an arbitrary skew base, the latter is the stronger concept.
\end{remark}

\begin{example}
A skew $\mathcal{V}$-category structure on a discrete category $\mathcal{A}$ is precisely a $\mathcal{V}$-category (in the usual sense) with set of objects $\mathcal{A}$. 
\end{example}

\begin{example} \label{closedex}
If $\mathcal{V}$ is a skew-monoidal category for which each functor $- \otimes X \colon \mathcal{V} \lra \mathcal{V}$ has a right adjoint $[X,-] \colon \mathcal{V}\lra \mathcal{V}$ (that is, $\mathcal{V}$ is a skew monoidal closed category \cite{MR3010098}), then $\mathcal{V}$ is skew-enriched over itself, with hom-objects $\underline{\mathcal{V}}(X,Y) = [X,Y]$. This skew $\mathcal{V}$-category $\underline{\mathcal{V}}$ is normal if and only if the skew-monoidal category $\mathcal{V}$ is left normal. 
\end{example}

\begin{remark} \label{weakmap}
Inspired by the examples of skew monoidal closed categories $\mathcal{V}$ in \cite{Bourke2015}, for which $\mathcal{V}(X,Y)$ and $\mathcal{V}(I \otimes X,Y) \cong \mathcal{V}(I,[X,Y])$ are respectively the sets of strict and weak morphisms from $X$ to $Y$ in $\mathcal{V}$, we can think of the elements of the sets  $\mathcal{A}(A,B)$ and $\mathcal{V}(I,\underline{\mathcal{A}}(A,B))$ respectively as the ``strict morphisms'' and ``weak morphisms'' from $A$ to $B$ in $\underline{\mathcal{A}}$. The functions (\ref{strwk}) then assign to each strict morphism a corresponding weak morphism.
\end{remark}

\begin{example} \label{skewexample}
Each small category $\mathcal{C}$ with set of objects $O$ determines a skew-monoidal structure on the functor category $[O,\Set]$, with tensor product $(X\otimes Y)_j = \sum_{i \in O} \mathcal{C}(i,j) \times X_i \times Y_j$ and unit object $1_j = 1$ \cite{MR3370862}. Each presheaf $F$ on $\mathcal{C}$ determines a skew-enrichment of any locally small category $\mathscr{E}$ over this skew-monoidal category, with hom-objects $\underline{\mathscr{E}}(A,B)_j = \mathscr{E}(A,B)^{Fj}$, composition
\begin{equation*}
\sum_i \mathcal{C}(i,j) \times \mathscr{E}(B,C)^{Fi} \times \mathscr{E}(A,B)^{Fj} \lra \mathscr{E}(A,C)^{Fj}
\end{equation*}
given by $(u,g,f) \mapsto (g_{(Fu)x} \circ f_x)_{x \in Fj}$, and unit $1 \lra \mathscr{E}(A,A)^{Fi}$ the constant at $1_A$.
\end{example}

\begin{example}
For any object $K$ of a locally small category $\mathscr{A}$ with powers and copowers, there is a skew monoidal closed structure on $\mathscr{A}$ with tensor product $A \otimes B = \mathscr{A}(K,B) \times A$, unit object $K$, and internal hom $[B,C] = C^{\mathscr{A}(K,B)}$. This left skew-monoidal structure is the reverse of the right skew-monoidal structure induced by the skew right warping $(\mathscr{A}(K,-),K,v,v_0,k)$ riding the action of $\Set$ on $\mathscr{A}$ by copowers \cite{MR3370862}, where  
\begin{equation*}
v_{A,B} \colon \mathscr{A}(K,A) \times \mathscr{A}(K,B) \lra \mathscr{A}(K,\mathscr{A}(K,A)\times B)
\end{equation*}
is the function that sends $(f,g)$ to the composite $\mathrm{in}_f  \circ g$, where $v_0 \colon 1 \lra \mathscr{A}(K,K)$ picks out $1_K$, and where $k_A \colon \mathscr{A}(K,A) \times K \lra A$ is the counit of the adjunction $-\times K \dashv \mathscr{A}(K,-)$.

For example, the category $\Set^{\mathcal{T}}$ of algebras in $\Set$ for a Lawvere theory $\mathcal{T}$ (e.g.\ the Lawvere theory for monoids, groups, abelian groups, rings, pointed sets, $M$-sets, etc.)\ is skew monoidal closed by taking $K$ to be the free $\mathcal{T}$-algebra on one generator; the tensor product classifies binary functions that are $\mathcal{T}$-homomorphisms in the first variable only, and the internal hom $[A,B]$ is the set of functions $\Set(UA,UB)$ between the underlying sets equipped with the pointwise $\mathcal{T}$-algebra structure. Moreover, for any locally small category $\mathscr{C}$ with finite products, the category of $\mathcal{T}$-algebras in $\mathscr{C}$ is skew-enriched over this skew-monoidal structure on $\Set^{\mathcal{T}}$, with enriched hom-objects $\underline{\mathscr{C}^\mathcal{T}}(A,B)$ the hom-sets $\mathscr{C}(A,B)$ equipped with the pointwise $\mathcal{T}$-algebra structure.
\end{example}

\begin{example} The $\mathscr{F}$-categories of \cite{MR2854177} are precisely the skew $\Cat$-categories for which each of the functions (\ref{strwk}) is injective. The tight and loose  morphisms of an $\mathscr{F}$-category are the strict and weak morphisms respectively of the corresponding skew $\Cat$-category. \end{example}

\begin{definition} A \emph{$\mathcal{V}$-functor} $(F,\underline{F}) \colon \underline{\mathcal{A}} \lra \underline{\mathcal{B}}$ between skew $\mathcal{V}$-categories consists of the following data:
\begin{enumerate}[(i)]
\item a functor $F \colon \mathcal{A} \lra  \mathcal{B}$,
\item a natural transformation $\underline{F} \colon \underline{\mathcal{A}}(A,B) \lra \underline{\mathcal{B}}(FA,FB)$,
\end{enumerate}
subject to the following two axioms.
\begin{equation*}
\cd{
\underline{\mathcal{A}}(B,C) \otimes \underline{\mathcal{A}}(A,B) \ar[d]_-{\underline{F} \otimes \underline{F}} \ar[r]^-M & \underline{\mathcal{A}}(A,C)  \ar[d]^-{\underline{F}} \\
\underline{\mathcal{B}}(FB,FC)\otimes \underline{\mathcal{B}}(FA,FB)  \ar[r]_-{M} & \underline{\mathcal{B}}(FA,FC)
}
\qquad
\cd{
I \ar[r]^-j \ar[dr]_-j & \underline{\mathcal{A}}(A,A) \ar[d]^-{\underline{F}} \\
& \underline{\mathcal{B}}(FA,FA)
}
\end{equation*}
\end{definition}

\begin{definition}
A \emph{$\mathcal{V}$-natural transformation} $\theta \colon (F,\underline{F}) \lra (G,\underline{G})$ between $\mathcal{V}$-functors consists of the following data:
\begin{enumerate}[(i)]
\item a natural transformation $\theta \colon F \lra G$,
\end{enumerate}
subject to the following axiom.
\begin{equation*}
\cd[@C=3em]{
\underline{\mathcal{A}}(A,B) \ar[r]^-{\underline{F}} \ar[d]_-{\underline{G}} & \underline{\mathcal{B}}(FA,FB) \ar[d]^-{\underline{\mathcal{B}}(1,\theta_B)} \\
\underline{\mathcal{B}}(GA,GB) \ar[r]_-{\underline{\mathcal{B}}(\theta_A,1)} & \underline{\mathcal{B}}(FA,GB)
}
\end{equation*}
\end{definition}

Note that, in the terminology of Remark \ref{weakmap}, the components of a $\mathcal{V}$-natural transformation are strict morphisms in $\underline{\mathcal{B}}$.

Skew $\mathcal{V}$-categories, $\mathcal{V}$-functors, and $\mathcal{V}$-natural transformations  form a $2$-category $\vskcat$, with the obvious compositions. Moreover, we have the following result, which we record for later reference.

\begin{proposition} \label{2fun} The assignment that sends a skew-monoidal category $\mathcal{V}$ to the $2$-category $\vskcat$ extends to a $2$-functor from the $2$-category of skew-monoidal categories to the $2$-category of ``large'' $2$-categories.
\begin{proof} The proof follows as in \cite{MR0225841}, with the further requirement that
change of base along a monoidal functor $F \colon \mathcal{V} \lra \mathcal{W}$ between skew-monoidal categories sends a skew $\mathcal{V}$-enrichment of a category $\mathcal{A}$ to a skew $\mathcal{W}$-enrichment of that same category $\mathcal{A}$.
\end{proof}
\end{proposition}

While every skew $\mathcal{V}$-category has an underlying $\mathcal{V}$-category ``of weak morphisms'' (by forgetting the morphisms of $\mathcal{A}$, the functoriality of $\underline{\mathcal{A}}(-,-)$, and the naturality of $M$ and $j$), to construct a skew $\mathcal{V}$-category from a $\mathcal{V}$-category requires the base skew-monoidal category $\mathcal{V}$ to be left normal. In that case, normal skew $\mathcal{V}$-categories are equivalent to $\mathcal{V}$-categories in the ordinary sense. Moreover, Richard Garner observed that in that case skew $\mathcal{V}$-categories are equivalently enriched categories over a certain base, to be described presently, and that the former equivalence follows from the latter by change of base.

For each skew-monoidal category $\mathcal{V}$, the functor $\mathcal{V}(I,-) \colon \mathcal{V} \lra \Set$ is monoidal, and hence there exists a unique skew-monoidal structure on the comma category $\mathscr{S}(\mathcal{V}) = \Set\downarrow\mathcal{V}(I,-)$ making the projection functors $P_1 \colon \mathscr{S}(\mathcal{V}) \lra \Set$ and $P_2 \colon \mathscr{S}(\mathcal{V}) \lra \mathcal{V}$ strict monoidal \cite{MR2167789}, which is left normal if $\mathcal{V}$ is so.

For each left normal skew-monoidal category $\mathcal{V}$, let $\vcat$ denote the $2$-category of $\mathcal{V}$-categories defined in \cite{MR3010098}.

\begin{proposition}
Let $\mathcal{V}$ be a left normal skew-monoidal category. Then the $2$-category $\vskcat$ of skew $\mathcal{V}$-categories is isomorphic to the $2$-category $\mathscr{S}(\mathcal{V})\text{-}\mathbf{Cat}$ of $\mathscr{S}(\mathcal{V})$-categories.
\begin{proof}
We first describe the 2-category $\mathscr{S}(\mathcal{V})\text{-}\mathbf{Cat}$ in elementary terms. An $\mathscr{S}(\mathcal{V})$-category amounts to a category $\mathcal{A}$ and a $\mathcal{V}$-category $\underline{\mathcal{A}}$ with the same set of objects, and an identity on objects functor from $\mathcal{A}$ to the underlying category of $\underline{\mathcal{A}}$. Such a functor amounts to a family of functions $\omega_{A,B} \colon \mathcal{A}(A,B) \lra \mathcal{V}(I,\underline{\mathcal{A}}(A,B))$ such that:\ (i) for each object $A$ we have $\omega_{A,A}(1_A) = j_A$, and (ii) for each composable pair of morphisms $f \colon A \lra B, g \colon B \lra C$ in the category $\mathcal{A}$, we have that $\omega(gf)$ is equal to the composite
\begin{equation*}
\cd[@C=3.7em]{
I \ar[r]^-r & I \otimes I \ar[r]^-{\omega(g) \otimes \omega(f)} & \underline{\mathcal{A}}(B,C) \otimes \underline{\mathcal{A}}(A,B) \ar[r]^-M & \underline{\mathcal{A}}(A,C).
}
\end{equation*}
An $\mathscr{S}(\mathcal{V})$-functor $(F,\underline{F}) \colon (\mathcal{A},\underline{\mathcal{A}},\omega) \lra (\mathcal{B},\underline{\mathcal{B}},\omega)$ amounts to a functor $F \colon \mathcal{A} \lra \mathcal{B}$ and a $\mathcal{V}$-functor  $\underline{F} \colon \underline{\mathcal{A}} \lra \underline{\mathcal{B}}$ (in the sense of \cite{MR3010098}) which agree on objects, such that $\omega(Ff) = \underline{F}\omega(f)$ for every morphism $f$ in $\mathcal{A}$. An $\mathscr{S}(\mathcal{V})$-natural transformation $\theta \colon (F,\underline{F}) \lra (G,\underline{G})$ amounts to a natural transformation $\theta \colon F \lra G$ such that $\omega(\theta) \colon \underline{F} \lra \underline{G}$ is $\mathcal{V}$-natural (in the sense of \cite{MR3010098}).

We now describe the isomorphism of $2$-categories. In one direction, a skew $\mathcal{V}$-category $(\mathcal{A},\underline{\mathcal{A}}(-,-))$ is sent to the $\mathscr{S}(\mathcal{V})$-category consisting of the category $\mathcal{A}$, the $\mathcal{V}$-category of weak morphisms in $\underline{\mathcal{A}}$ described above, and the family of functions (\ref{strwk}). Similarly, $\mathcal{V}$-functors and $\mathcal{V}$-natural transformations are sent to their evident underlying $\mathscr{S}(\mathcal{V})$-functors and $\mathscr{S}(\mathcal{V})$-natural transformations.

In the other direction, an $\mathscr{S}(\mathcal{V})$-category $(\mathcal{A},\underline{\mathcal{A}},\omega)$ is sent to a skew $\mathcal{V}$-enrichment of the category $\mathcal{A}$ with enriched hom-objects $\underline{\mathcal{A}}(A,B)$. Using left normality of $\mathcal{V}$, we can extend $\underline{\mathcal{A}}(-,-)$ to a functor $\mathcal{A}^{\mathrm{op}} \times \mathcal{A} \lra \mathcal{V}$ making $M$ and $j$ natural. For each $g \colon B \lra C$ in $\mathcal{A}$, we define $\underline{\mathcal{A}}(A,g)$ to be the composite
\begin{equation*}
\cd[@C=3em]{
\underline{\mathcal{A}}(A,B) \ar[r]^-{l^{-1}} & I \otimes \underline{\mathcal{A}}(A,B) \ar[r]^-{\omega(g) \otimes 1} & \underline{\mathcal{A}}(B,C) \otimes \underline{\mathcal{A}}(A,B) \ar[r]^-M & \underline{\mathcal{A}}(A,C),
}
\end{equation*}
and for each $f \colon A \lra B$ in $\mathcal{A}$, we define $\underline{\mathcal{A}}(f,C)$ to be the composite 
\begin{equation*}
\cd[@C=3em]{
\underline{\mathcal{A}}(B,C) \ar[r]^-r & \underline{\mathcal{A}}(B,C) \otimes I \ar[r]^-{1 \otimes \omega(f)} & \underline{\mathcal{A}}(B,C) \otimes \underline{\mathcal{A}}(A,B) \ar[r]^-M & \underline{\mathcal{A}}(A,C).
}
\end{equation*}
Moreover, the components of $\mathscr{S}(\mathcal{V})$-functors and $\mathscr{S}(\mathcal{V})$-natural transformations are natural with respect to the morphisms of $\mathcal{A}$ under this definition of the functor $\underline{\mathcal{A}}(-,-)$.

Both assignments are easily seen to define mutually inverse $2$-functors.
\end{proof}
\end{proposition}

\begin{corollary} \label{2adj}
Let $\mathcal{V}$ be a left normal skew-monoidal category. Then there is a $2$-adjunction
\begin{equation*}
\cd{
\vcat \ar@<-1.5ex>[rr]^-{\hdash} && \ar@<-1.5ex>[ll] \vskcat
}
\end{equation*}
whose counit is an identity, and whose right $2$-adjoint restricts to a $2$-equivalence between $\vcat$ and the full sub-$2$-category of $\vskcat$ on the normal skew $\mathcal{V}$-categories.
\begin{proof}
Identifying $\vskcat$ with the isomorphic $2$-category $\mathscr{S}(\mathcal{V})\text{-}\mathbf{Cat}$, this 2-adjunction is the image under the change of base $2$-functor from Proposition \ref{2fun} of the monoidal adjunction 
\begin{equation*}
\cd{
\mathcal{V} \ar@<-1.5ex>[rr]^-{\hdash} && \ar@<-1.5ex>[ll] \mathscr{S}(\mathcal{V})
}
\end{equation*}
whose left adjoint is the strict monoidal projection functor $P_2$, and whose right adjoint sends an object $X$ of $\mathcal{V}$ to the object $(\mathcal{V}(I,X),1_{\mathcal{V}(I,X)},X)$ of the comma category $\mathscr{S}(\mathcal{V})$. The counit of this monoidal adjunction is an identity and hence so is the counit of the induced $2$-adjunction.

The component of the unit of this monoidal adjunction at an object $(E,f,X)$ of $\mathscr{S}(\mathcal{V})$ is an isomorphism if and only if the function $f \colon E \lra \mathcal{V}(I,X)$ is a bijection. Hence the component of the unit of the induced $2$-adjunction at an $\mathscr{S}(\mathcal{V})$-category $(\mathcal{A},\underline{\mathcal{A}},\omega)$ is an isomorphism if and only if all the functions $\omega \colon \mathcal{A}(A,B) \lra \mathcal{V}(I,\underline{\mathcal{A}}(A,B))$ are bijections, that is, if the corresponding skew $\mathcal{V}$-category is normal.
\end{proof}
\end{corollary}

\section{Locally weak comonads} \label{locwksection}

In this section we define the notion of locally weak comonad, which is fundamental to the theory of enriched algebraic weak factorisation systems, and show that a locally weak comonad on a skew-enriched category gives rise to a new skew-enrichment  of the same category. This result is the analogue for skew-enriched categories of a result for monoidal comonads on skew-monoidal categories \cite{MR2964621}  and a result for closed comonads on skew-closed categories \cite{MR3010098}, which we combine as Proposition \ref{skmonthm} below.

\begin{proposition}[Szlach{\'a}nyi, Street] \label{skmonthm}
Let $(Q,\varphi,\varphi_0,\delta,\varepsilon)$ be a monoidal comonad on a skew-monoidal category $\mathcal{V}$. Then the following data define a skew-monoidal structure on $\mathcal{V}$:
\begin{enumerate}[font=\normalfont, label=(\roman*)]
\item tensor product functor $X \otimes QY$, 
\item unit object $I$, 
\item associativity natural transformation
\begin{equation*}
\cd[@C=2.75em]{
(X \otimes QY) \otimes QZ \ar[r]^-a & X \otimes (QY \otimes QZ) \ar[r]^-{1\otimes (1 \otimes \delta)} & X \otimes (QY \otimes Q^2Z) \ar[r]^-{1 \otimes \varphi} & X \otimes Q(Y\otimes QZ),
}
\end{equation*}
\item left unit natural transformation
\begin{equation*}
\cd{
I \otimes QX \ar[r]^-l & QX \ar[r]^-{\varepsilon} & X,
}
\end{equation*}
\item right unit natural transformation
\begin{equation*}
\cd{
X \ar[r]^--r & X \otimes I \ar[r]^-{1 \otimes \varphi_0} & X \otimes QI.
}
\end{equation*}
\end{enumerate}
Moreover, if the original skew-monoidal structure on $\mathcal{V}$ is closed with internal hom functor $[X,Y]$, then the skew-monoidal structure defined above is closed with internal hom functor $[QX,Y]$.
\end{proposition}

\begin{example} 
The strictification comonad $\textbf{st}$  on the category $\twocat$ of 2-categories and 2-functors, which sends a 2-category to its strictification (a.k.a.\ its pseudofunctor classifier), extends to a monoidal closed comonad with respect to the Gray symmetric monoidal closed structure on $\twocat$ \cite{MR1261589}. By the previous proposition, this induces a skew monoidal closed structure on $\twocat$ (considered in \cite{Bourke2015}), whose internal hom $\textbf{Gray}(\textbf{st}A,B)$ is the 2-category of pseudofunctors $A \lra B$, pseudonatural transformations, and modifications.
\end{example}

\begin{remark} \label{weakmaps2}
More generally, given a monoidal comonad $Q$ on a skew-monoidal category $\mathcal{V}$, we could think of the co-Kleisli morphisms of $Q$ as ``weak morphisms'' in $\mathcal{V}$; in the previous example, the weak morphisms are the pseudofunctors. By Proposition \ref{skmonthm}, this agrees with our previous notion of weak morphism (Remark \ref{weakmap}) at least when the original skew-monoidal structure on $\mathcal{V}$ is left normal. This informs our choice of terminology in the following definition.
\end{remark}

Let $(\mathcal{V},\otimes,I)$ be a skew-monoidal category.
\begin{definition}
A \emph{locally weak $\mathcal{V}$-comonad} on a skew $\mathcal{V}$-category $\underline{\mathcal{A}}$ consists of the following data:
\begin{enumerate}[(i)]
\item a monoidal comonad $(Q,\varphi,\varphi_0,\delta,\varepsilon)$ on $\mathcal{V}$,
\item a comonad $(S,\delta,\varepsilon)$ on the category $\mathcal{A}$,
\item a natural transformation $\psi \colon Q\underline{\mathcal{A}}(A,B) \lra \underline{\mathcal{A}}(SA,SB)$,
\end{enumerate}
subject to the following four axioms.
\begin{equation*}
\cd{
Q\underline{\mathcal{A}}(B,C) \otimes Q\underline{\mathcal{A}}(A,B) \ar[r]^-{\varphi} \ar[d]_-{\psi \otimes \psi} & Q(\underline{\mathcal{A}}(B,C) \otimes \underline{\mathcal{A}}(A,B)) \ar[r]^-{QM} & Q\underline{\mathcal{A}}(A,C) \ar[d]^-{\psi} \\
\underline{\mathcal{A}}(SB,SC) \otimes \underline{\mathcal{A}}(SA,SB) \ar[rr]_-M && \underline{\mathcal{A}}(SA,SC)
}
\end{equation*}
\begin{equation*}
\cd{
Q^2\underline{\mathcal{A}}(A,B) \ar[r]^-{Q\psi} & Q\underline{\mathcal{A}}(SA,SB) \ar[r]^-{\psi} & \underline{\mathcal{A}}(S^2A,S^2B) \ar[d]^-{\underline{\mathcal{A}}(\delta,1)} \\
Q\underline{\mathcal{A}}(A,B) \ar[u]^-{\delta} \ar[r]_-{\psi} & \underline{\mathcal{A}}(SA,SB) \ar[r]_-{\underline{\mathcal{A}}(1,\delta)} & \underline{\mathcal{A}}(SA,S^2B)
}
\end{equation*}
\begin{equation*}
\cd{
QI \ar[r]^-{Qj} & Q\underline{\mathcal{A}}(A,A) \ar[d]^-{\psi} \\
I \ar[u]^-{\varphi_0} \ar[r]_-j & \underline{\mathcal{A}}(SA,SA)
}
\qquad
\cd{
Q\underline{\mathcal{A}}(A,B) \ar[r]^-{\psi} \ar[d]_{\varepsilon} & \underline{\mathcal{A}}(SA,SB) \ar[d]^-{\underline{\mathcal{A}}(1,\varepsilon)} \\
\underline{\mathcal{A}}(A,B) \ar[r]_-{\underline{\mathcal{A}}(\varepsilon,1)} & \underline{\mathcal{A}}(SA,B)
}
\end{equation*}
\end{definition}

We will sometimes call a locally weak $\mathcal{V}$-comonad a \emph{locally $Q$-weak $\mathcal{V}$-comonad} to indicate its dependence on its underlying monoidal comonad $Q$.

\begin{remark}
Note that when $Q$ is the identity monoidal comonad on $\mathcal{V}$, a locally $Q$-weak $\mathcal{V}$-comonad is simply a comonad in the 2-category $\vskcat$. In the general case, $S$ fails to be a $\mathcal{V}$-functor because its action on homs $\psi \colon Q\underline{\mathcal{A}}(A,B) \lra \underline{\mathcal{A}}(SA,SB)$ is given by weak morphisms (in the sense of Remark \ref{weakmaps2}). Hence the name ``locally weak $\mathcal{V}$-comonad''.
\end{remark}

\begin{example}
Any monoidal closed comonad $Q$ on a skew monoidal closed category $\mathcal{V}$ is a locally $Q$-weak $\mathcal{V}$-comonad on the skew $\mathcal{V}$-category $\underline{\mathcal{V}}$. (Note that any monoidal comonad or closed comonad on such a $\mathcal{V}$ can be uniquely extended to a monoidal closed comonad on $\mathcal{V}$ by adjointness.)
In particular, the strictification comonad $\textbf{st}$ on $\twocat$ is a locally weak $\textbf{Gray}$-comonad on the $\textbf{Gray}$-category $\textbf{Gray}$ (the self-enrichment of $\twocat$ over the Gray symmetric monoidal closed structure). By the universal property of $\textbf{st}$ as the pseudofunctor classifier comonad, the locally weak $\textbf{Gray}$-functor structure morphisms for $\textbf{st}$ correspond to pseudofunctors $\textbf{Gray}(A,B) \lra \textbf{Gray}(\textbf{st}A,\textbf{st}B)$.
\end{example}

The notion of locally weak comonad enjoys the following two abstract descriptions (the second of which was pointed out to the author by Richard Garner). A locally weak $\mathcal{V}$-comonad (with underlying monoidal comonad $Q$) is equally:
\begin{enumerate}[leftmargin=*, label=(\roman*)]
\item a comonad in the 2-category $\skenrcat$ of skew-enriched categories over an arbitrary skew-monoidal base, where $\skenrcat \lra \skmoncat$ is the $2$-opfibration corresponding to the change of base 2-functor $\skmoncat \lra \bigtwocat$ of Proposition \ref{2fun} (see \cite{MR3153612});
\item a comonad in the 2-category of skew $\mathcal{V}$-categories and ``locally weak $\mathcal{V}$-functors'', that is, the co-Kleisli 2-category for the 2-comonad $Q_{\ast}$ on $\vcat$ defined by change of base along $Q$ (which is the image of the monoidal comonad $Q$ under the change of base 2-functor  $\skmoncat \lra \bigtwocat$).
\end{enumerate}

\begin{proposition} \label{locwkthm}
Let $(Q,S,\psi)$ be a locally weak $\mathcal{V}$-comonad on a skew $\mathcal{V}$-category \linebreak $(\mathcal{A},\underline{\mathcal{A}}(-,-),M,j)$. Then the following data define a skew-enrichment of the category $\mathcal{A}$ over the skew-monoidal structure on $\mathcal{V}$ induced by the monoidal comonad $Q$  (described in Proposition \ref{skmonthm}):
\begin{enumerate}[{\normalfont (i)}]
\item hom functor $\underline{\mathcal{A}}(SA,B)$,
\item composition natural transformation
\fontsize{11pt}{0}
\begin{equation*}
\cd[@C=3em]{
\underline{\mathcal{A}}(SB,C) \otimes Q\underline{\mathcal{A}}(SA,B) \ar[r]^-{1 \otimes \psi} & \underline{\mathcal{A}}(SB,C)\otimes \underline{\mathcal{A}}(S^2A,SB) \ar[r]^-M & \underline{\mathcal{A}}(S^2A,C) \ar[r]^-{\underline{\mathcal{A}}(\delta,1)} & \underline{\mathcal{A}}(SA,C),
}
\end{equation*}
\normalsize
\item unit natural transformation
\begin{equation*}
\cd[@C=3em]{
I \ar[r]^-j & \underline{\mathcal{A}}(A,A) \ar[r]^-{\underline{\mathcal{A}}(\varepsilon,1)} & \underline{\mathcal{A}}(SA,A). 
}
\end{equation*}
\end{enumerate}
\begin{proof}
The diagram for associativity is displayed on the following page, where we write $X \cdot Y$ for $X\otimes Y$ and $(A,B)$ for $\underline{\mathcal{A}}(A,B)$.

\afterpage{
\begin{landscape}
\fontsize{7.5pt}{0}
\vspace*{\fill}
\begin{center}
\begin{equation*}
\cd[@C=2em]{
(SC,D)\cdot (Q(SB,C)\cdot Q(SA,B))  \ar[rrr]^-{1 \cdot (1\cdot\delta)} &&& (SC,D)\cdot(Q(SB,C)\cdot Q^2(SA,B)) \ar[r]^-{1\cdot \varphi} \ar[d]^-{1\cdot(1\cdot Q\psi)} & (SC,D) \cdot Q((SB,C) \cdot Q(SA,B)) \ar[d]^-{1\cdot Q(1\cdot \psi)}  \\
((SC,D) \cdot Q(SB,C))\cdot Q(SA,B) \ar[u]^-a  \ar[r]^-{1\cdot \delta} \ar[d]_-{1 \cdot \psi} & ((SC,D)\cdot Q(SB,C))\cdot Q^2(SA,B) \ar[r]^-{1 \cdot Q\psi} & ((SC,D)\cdot Q(SB,C)) \cdot Q(S^2A,SB) \ar[r]^-a \ar[d]^-{1\cdot \psi} & (SC,D)\cdot (Q(SB,C)\cdot Q(S^2A,SB)) \ar[r]^-{1 \cdot \psi} \ar[dd]^-{1 \cdot (\psi \cdot \psi)}  & (SC,D)\cdot Q((SB,C)\cdot (S^2A,SB)) \ar[dddd]^-{1 \cdot QM}  \\
((SC,D)\cdot Q(SB,C)) \cdot (S^2A,SB) \ar[r]^-{1 \cdot (1,\delta)} \ar[d]_-{(1\cdot \psi)\cdot 1} & ((SC,D) \cdot Q(SB,C)) \cdot (S^2A,S^2B) \ar[d]_-{(1\cdot \psi)\cdot 1} 
 & ((SC,D) \cdot Q(SB,C)) \cdot (S^3A,S^2B) \ar[l]_-{1\cdot (\delta,1)} \ar[d]^-{(1\cdot \psi) \cdot 1} \\
((SC,D) \cdot (S^2B,SC)) \cdot (S^2A,SB) \ar[r]^-{1\cdot (1,\delta)} \ar[d]_-{M\cdot 1} & ((SC,D)\cdot (S^2B,SC)) \cdot (S^2A,S^2B) \ar[d]_-{M\cdot 1} & ((SC,D) \cdot (S^2B,SC)) \cdot (S^3A,S^2B) \ar[l]_-{1 \cdot (\delta,1)} \ar[d]^-{M \cdot 1} \ar[r]^-a & (SC,D) \cdot ((S^2B,SC) \cdot (S^3A,S^2B)) \ar[dd]^-{1 \cdot M} \\
(S^2B,D) \cdot (S^2A,SB) \ar[r]^-{1 \cdot (1,\delta)} \ar[d]_-{(\delta,1) \cdot 1} & (S^2B,D) \cdot (S^2A,S^2B)  \ar[d]_-M & (S^2B,D) \cdot (S^3A,S^2B)  \ar[l]_-{1 \cdot (\delta,1)} \ar[d]^-M \\
(SB,D) \cdot (S^2A,SB) \ar[r]_-M & (S^2A,D) \ar[d]_-{(\delta,1)}& (S^3A,D) \ar[d]^-{(S\delta,1)} \ar[l]^-{(\delta,1)} & (SC,D) \cdot (S^3A,SC)  \ar[l]^-M \ar[d]^-{1\cdot (S\delta,1)}& (SC,D) \cdot Q(S^2A,C) \ar[l]^-{1\cdot \psi} \ar[d]^-{1\cdot Q(\delta,1)} \\
& (SA,D) & (S^2A,D) \ar[l]^-{(\delta,1)} & (SC,D) \cdot (S^2A,SC) \ar[l]^-M & (SC,D) \cdot Q(SA,C) \ar[l]^-{1\cdot \psi}
}
\end{equation*}
\end{center}
\vspace*{\fill}
\end{landscape}
}

\begin{equation*}
\cd[@C=3em]{
I \otimes Q\underline{\mathcal{A}}(SA,B) \ar[r]^-{1 \otimes \varepsilon} \ar[d]_-{1 \otimes \psi} & I \otimes \underline{\mathcal{A}}(SA,B) \ar[r]^-l \ar[d]_-{1 \otimes \underline{\mathcal{A}}(\varepsilon,1)} & \underline{\mathcal{A}}(SA,B) \ar[d]_-{\underline{\mathcal{A}}(\varepsilon,1)} \ar@/^3pc/[ddd]^-1 \\
I \otimes \underline{\mathcal{A}}(S^2A,SB) \ar[r]^-{1\otimes \underline{\mathcal{A}}(1,\varepsilon)} \ar[d]_-{j \otimes 1} & I \otimes \underline{\mathcal{A}}(S^2A,B) \ar[r]^-l \ar[d]_-{j \otimes 1} & \underline{\mathcal{A}}(S^2A,B) \ar[d]_-1 \\
\underline{\mathcal{A}}(B,B) \otimes \underline{\mathcal{A}}(S^2A,SB) \ar[r]^-{1\otimes \underline{\mathcal{A}}(1,\varepsilon)} \ar[d]_-{\underline{\mathcal{A}}(\varepsilon,1)\otimes 1} & \underline{\mathcal{A}}(B,B) \otimes \underline{\mathcal{A}}(S^2A,B) \ar[r]^-M \ar[d]_-M & \underline{\mathcal{A}}(S^2A,B) \ar[d]_-{\underline{\mathcal{A}}(\delta,1)} \\
\underline{\mathcal{A}}(SB,B) \otimes \underline{\mathcal{A}}(S^2A,SB) \ar[r]_-M & \underline{\mathcal{A}}(S^2A,B) \ar[r]_-{\underline{\mathcal{A}}(\delta,1)} & \underline{\mathcal{A}}(SA,B)
}
\end{equation*}

\normalsize

\begin{equation*}
\begin{gathered}[b]
\cd[@C=3.75em]{
\underline{\mathcal{A}}(SA,B) \otimes QI \ar[r]^-{1\otimes Qj} & \underline{\mathcal{A}}(SA,B) \otimes Q\underline{\mathcal{A}}(A,A) \ar[r]^-{1 \otimes Q\underline{\mathcal{A}}(\varepsilon,1)} \ar[d]^-{1 \otimes \psi} & \underline{\mathcal{A}}(SA,B) \otimes Q\underline{\mathcal{A}}(SA,A) \ar[d]^-{1 \otimes\psi} \\
\underline{\mathcal{A}}(SA,B) \otimes I \ar[u]^-{1 \otimes \varphi_0} \ar[r]^-{1 \otimes j} & \underline{\mathcal{A}}(SA,B) \otimes \underline{\mathcal{A}}(SA,SA) \ar[r]^-{1 \otimes \underline{\mathcal{A}}(S\varepsilon,1)} \ar[d]^-M & \underline{\mathcal{A}}(SA,B) \otimes \underline{\mathcal{A}}(S^2A,SA) \ar[d]^-M \\
\underline{\mathcal{A}}(SA,B) \ar[u]^-r \ar[r]_-1 & \underline{\mathcal{A}}(SA,B) \ar[r]_-{\underline{\mathcal{A}}(S\varepsilon,1)} \ar[dr]_-1 & \underline{\mathcal{A}}(S^2A,B) \ar[d]^-{\underline{\mathcal{A}}(\delta,1)} \\
&& \underline{\mathcal{A}}(SA,B) 
}\\[-\dp\strutbox]
\end{gathered}
\qedhere
\end{equation*} 
\end{proof}
\end{proposition}

\pagebreak 
\begin{example}[See \cite{ct2017talk}]
Let $(L,R)$ be a monoidal algebraic weak factorisation system on a (suitably nice) monoidal category $\mathcal{V}$. Then the  cofibrant replacement comonad $Q$ for $(L,R)$ is a monoidal comonad on $\mathcal{V}$. Furthermore, the cofibrant replacement comonad $S$ for an $(L,R)$-enriched algebraic weak factorisation system $(H,M)$ on a $\mathcal{V}$-category $\mathcal{A}$ is a locally $Q$-weak $\mathcal{V}$-comonad on $\mathcal{A}$. Proposition \ref{locwkthm} then defines a skew-enrichment of (the underlying category of) $\mathcal{A}$ over the skew-monoidal structure on $\mathcal{V}$ induced by $Q$, from which it follows that the co-Kleisli category for $S$ (known as the category of weak maps for $(H,M)$ \cite{MR3393454}) is enriched over this skew-monoidal structure on $\mathcal{V}$.
\end{example}

\section{Skew-proactegories} \label{proactsection}
In this section we prove that the structure of a skew $\mathcal{V}$-category $(\mathcal{A},\underline{\mathcal{A}}(-,-))$ for which every functor $\underline{\mathcal{A}}(A,-) \colon \mathcal{A} \lra \mathcal{V}$ is a right adjoint is equivalent to the structure of a skew $\mathcal{V}$-actegory \cite{MR3370862} $(\mathcal{A},\ast)$ for which every functor $-\ast A \colon \mathcal{V} \lra \mathcal{A}$ is a left adjoint. This is an extension to the skew setting of the ``often-rediscovered folklore'' equivalence, for $\mathcal{V}$ a monoidal closed category, between tensored $\mathcal{V}$-categories and hommed $\mathcal{V}$-actegories \cite{MR1466618,MR1897810}. 

Our method of proof is an adaptation of a folklore argument that proves the equivalence between closed skew-monoidal categories and monoidal skew-closed categories (itself an extension to the skew setting of a result of \cite{MR0225841}) as a consequence of the result that both skew-monoidal categories and skew-closed categories are equivalent to special kinds of skew-promonoidal category \cite{MR3010098}. Therefore we introduce the notion of skew $\mathcal{V}$-proactegory, and show that both skew $\mathcal{V}$-categories and skew $\mathcal{V}$-actegories are equivalent to special kinds of skew $\mathcal{V}$-proactegory.

Note that throughout this section we take $\Set$ to be the category of sets belonging to some universe with respect to which our categories bearing skew-promonoidal or skew-proactegory structures are small.

Let $(\mathcal{V},P,J,\alpha,\lambda,\rho)$ be a left skew-promonoidal category \cite{MR3010098}.

\begin{definition}
A \emph{left skew $\mathcal{V}$-proactegory} $(\mathcal{A},T,\alpha,\lambda)$ consists of the following data:
\begin{enumerate}[(i)]
\item a category $\mathcal{A}$,
\item a functor $T \colon \mathcal{V}^\text{op} \times \mathcal{A}^\mathrm{op} \times \mathcal{A} \lra \Set$,
\item a natural transformation $$\alpha \colon \int^B T(X,B;C) \times T(Y,A;B) \lra \int^Z T(Z,A;C) \times P(X,Y;Z),$$
\item a natural transformation
$$\lambda \colon \mathcal{A}(A,B) \lra \int^X T(X,A;B)\times JX,$$
\end{enumerate}
subject to the following three axioms.

\small
\begin{equation*}
\cd[@C=3.5em]{
\int^{B,C} T(X,C;D) \times T(Y,B;C) \times T(Z,A;B) \ar[r]^-{\int^B \alpha \times 1} \ar[d]_-{\int^C 1 \times \alpha} & \int^{U,B} T(U,B;D) \times P(X,Y;U) \times T(Z,A;B) \ar[d]^-{\cong} \\
\int^{V,C} T(X,C;D) \times T(V,A;C) \times P(Y,Z;V) \ar[d]_-{\int^V \alpha \times 1} & \int^{U,B} T(U,B;D) \times T(Z,A;B) \times P(X,Y;U) \ar[d]^-{\int^U \alpha \times 1} \\
\int^{V,W} T(W,A;D) \times P(X,V;W) \times P(Y,Z;V) \ar[r]_-{\int^W 1 \times \alpha} & \int^{U,W} T(W,A;D) \times P(U,Z;W) \times P(X,Y;U)
}
\end{equation*}

\normalsize
\begin{equation*}
\cd[@C=3.5em]{
\int^C \mathcal{A}(C,B) \times T(X,A;C) \ar[r]^-{\int^C \lambda \times 1} \ar[d]_-{\cong} & \int^{U,C} T(U,C;B) \times JU \times T(X,A;C) \ar[d]^-{\cong} \\
T(X,A;B) \ar[d]_-{\cong} & \int^{U,C} T(U,C;B) \times T(X,A;C) \times JU \ar[d]^-{\int^U \alpha \times 1} \\
\int^W T(W,A;B) \times \mathcal{V}(X,W) \ar[r]_-{\int^W 1 \times \lambda} & \int^{U,W}T(W,A;B) \times P(U,X;W) \times JU
}
\end{equation*}

\begin{equation*}
\cd[@C=3.5em]{
\int^C T(X,C;B) \times \mathcal{A}(A,C) \ar[r]^-{\cong} \ar[d]_-{\int^C 1 \times \lambda} & T(X,A;B) \ar[dd]^-{\cong} \\
\int^{V,C} T(X,C;B) \times T(V,A;C) \times JV \ar[d]_-{\int^V \alpha \times 1} \\
\int^{V,W} T(W,A;B) \times P(X,V;W) \times JV \ar[r]_-{\int^W 1 \times \rho} & \int^W T(W,A;B) \times \mathcal{V}(X,W)
}
\end{equation*}

A \emph{$\mathcal{V}$-proactegory} is a left skew $\mathcal{V}$-proactegory for which $\alpha$ and $\lambda$ are invertible.
\end{definition}

One can define similarly the notion of right skew $\mathcal{V}$-proactegory, for $\mathcal{V}$ a right skew-pro\-monoidal category, by reversing the directions of the natural transformations $\alpha$ and $\lambda$ in the above definition.

A left (right) skew $\mathcal{V}$-proactegory is precisely a right (left) skew-action of the right (left) skew-monoidale $\mathcal{V}$ \cite{MR2972969,RamonPhD} in the monoidal bicategory of categories and profunctors, under the convention that a profunctor from $\mathcal{A}$ to $\mathcal{B}$ is a functor $\mathcal{A}^\mathrm{op}\times \mathcal{B} \lra \Set$.

\begin{example}
Every skew-promonoidal category $(\mathcal{V},P)$ is canonically a skew $\mathcal{V}$-proactegory with skew-proaction $T=P$.
\end{example}

\begin{example} \label{skewexample2}
A right skew-promonoidal structure on a set $O$ for which $Ji = 1$ for all $i \in O$ amounts precisely to a small category with set of objects $O$. The right skew-promonoidal structure determined by such a category $\mathcal{C}$ has $P(i,j;k) = \mathcal{C}(i,j)$ if $j=k$ and empty otherwise, $Ji = 1$ for all $i \in O$, associativity constraint $\alpha$ given by the functions $\mathcal{C}(j,k) \times \mathcal{C}(i,j) \lra \mathcal{C}(i,k) \times \mathcal{C}(j,k)$ sending $(v,u)$ to $(vu,v)$, left unit constraint $\lambda$ given by the unique functions $\sum_i \mathcal{C}(i,j) \lra 1$, and right unit constraint $\rho$ given by the functions $1 \lra \mathcal{C}(i,i)$ picking out the identity $1_i$. Moreover, the category $\mathcal{C}$ is a groupoid if and only if the associativity constraint $\alpha$ is invertible.

A right skew-proaction on the terminal category $1$ of the right skew-promonoidal category determined as above by a small category $\mathcal{C}$ amounts precisely to a presheaf $F$ on $\mathcal{C}$, with $T(i,\ast;\ast) = Fi$, associativity constraint $\alpha$ given by the functions $Fj \times \mathcal{C}(i,j) \lra Fi \times Fj$ sending $(x,u)$ to $((Fu)x,x)$, and unit constraint $\lambda$ given by the unique function $\sum_i Fi \lra 1$. Moreover, if $\mathcal{C}$ is a groupoid, then the presheaf $F$ is a torsor if and only if the associativity constraint $\alpha$ is invertible and the unit constraint $\lambda$ is surjective.
\end{example}

Before commencing the main argument of this section, we give a brief treatment of the analogues of Day convolution \cite{DayPhD} for skew-proactegories. These take the form of ``convolution skew-actions'' (Proposition \ref{dayconvprop}) and ``convolution skew-enrichments'' (Proposition \ref{convenr}).

By Day convolution, to give a right (left) skew-pro\-monoidal structure on a category $\mathcal{V}$ is equivalently to give a left (right) skew-monoidal structure on the functor category $[\mathcal{V},\Set]$ whose tensor product functor preserves colimits in each variable. The proof of this statement can be adapted to prove the following result.

\begin{proposition} \label{dayconvprop}
Let $\mathcal{V}$ be a right (left) skew-promonoidal category, and let $(\mathcal{A},T)$ be a right (left) skew $\mathcal{V}$-proactegory. Then the functor category $[\mathcal{A},\Set]$ is a left (right) skew $[\mathcal{V},\Set]$-actegory whose action functor, which is given by the formula
\begin{equation*}
M \ast F = \int^{X,A} T(X,A;-) \times MX \times FA,
\end{equation*}
 preserves colimits in each variable.
 Moreover, every such left (right) skew $[\mathcal{V},\Set]$-actegory structure on $[\mathcal{A},\Set]$ arises in this way from a right (left) skew $\mathcal{V}$-proactegory structure on $\mathcal{A}$.
\end{proposition}

More generally, if $\mathscr{X}$ is a cocomplete left (right) skew-monoidal category whose tensor product functor preserves colimits in each variable, and $\mathscr{C}$ is a  cocomplete left (right) skew $\mathscr{X}$-actegory, then similar formulas define a left (right) skew action of $[\mathcal{V},\mathscr{X}]$, equipped with the Day convolution left (right) skew-monoidal structure, on the functor category $[\mathcal{A},\mathscr{C}]$, such that the action functor preserves colimits in each variable. We state our result on convolution skew-enrichments at this level of generality.

\begin{proposition} \label{convenr}
Let $\mathcal{V}$ be a right skew-promonoidal category, and let $(\mathcal{A},T)$ be a right skew $\mathcal{V}$-proactegory. Moreover, let $\mathscr{X}$ be a  complete and cocomplete left skew-monoidal category whose tensor product functor preserves colimits in each variable,  and let $(\mathscr{C},\underline{\mathscr{C}}(-,-))$ be a skew $\mathscr{X}$-category. Then the functor category $[\mathcal{A},\mathscr{C}]$ is skew-enriched over  $[\mathcal{V},\mathscr{X}]$, equipped with the Day convolution left skew-monoidal structure, with hom-objects given by the following formula (where $\pitchfork$ denotes powers in $\mathscr{X}$).
\begin{equation*}
\underline{\Hom}(F,G) = \int_{A,B} T(-,A;B) \pitchfork \underline{\mathscr{C}}(FA,GB)
\end{equation*}
\end{proposition}

\begin{example}
Let $\mathcal{C}$ be a small category with set of objects $O$. The skew-enrichment of a locally small category $\mathscr{E}$ over the skew-monoidal structure on $[O,\Set]$ determined by a presheaf on $\mathcal{C}$ (Example \ref{skewexample}) is the convolution skew-enrichment for the right skew-proaction of $O$ on $1$ determined by that presheaf (Example \ref{skewexample2}).
\end{example}

\begin{example}
In the ordinary (non-skew) setting, the special case of Proposition \ref{convenr} in which $\mathcal{A}$ is the canonical proaction of $\mathcal{V}$ on itself is a classical result (as are its enriched generalisations). Familiar examples  include:\ (1) the enrichment of the category of simplicial objects in any locally small category over the cartesian monoidal category of simplicial sets, (2) the enrichment of the category of chain complexes in any additive category over the symmetric monoidal closed additive category of chain complexes of abelian groups, and (3) for $\mathscr{V}$ a sufficiently nice monoidal category and $\mathscr{C}$ a $\mathscr{V}$-category,  the enrichment of the arrow category $[\mathbf{2},\mathscr{C}_0]$ over the arrow category $[\mathbf{2},\mathscr{V}]$, with tensor product and hom-objects given respectively by the Leibniz constructions of pushout-product and pullback-hom (see \cite{MR3217884,ct2017talk}).
\end{example}

We now commence the main argument of this section  by  showing that both skew-enriched categories and skew-actegories are equivalent to special kinds of skew-proactegory. Note that for the remainder of this section, all skew structures will be left skew.

\begin{proposition} \label{yonedaprop} Let $\mathcal{V}$ be a skew-monoidal category, let $\mathcal{A}$ be  a category, and let \linebreak $T \colon \mathcal{V}^\mathrm{op} \times \mathcal{A}^\mathrm{op} \times \mathcal{A} \lra \Set$ be a functor. Then the following are true.
\begin{enumerate}[leftmargin=*, font=\normalfont, label=(\alph*)]
\item Given a representation
\begin{equation} \label{actiso}
T(X,A;-) \cong \mathcal{A}(X \ast A,-) : \mathcal{A} \lra \Set
\end{equation}
 for each $X\in\mathcal{V}, A \in \mathcal{A}$, 
there is a bijection between skew $\mathcal{V}$-proactegory structures extending $(\mathcal{A},T)$ and skew $\mathcal{V}$-actegory structures extending $(\mathcal{A},\ast)$ such that the isomorphisms \emph{(\ref{actiso})} are natural in each variable.
\item Given a representation
\begin{equation} \label{homiso}
T(-,A;B) \cong \mathcal{V}(-,\underline{\mathcal{A}}(A,B)) : \mathcal{V}^\mathrm{op} \lra \Set
\end{equation}
 for each $A,B \in \mathcal{A}$,
there is a bijection between skew $\mathcal{V}$-proactegory structures extending $(\mathcal{A},T)$ and skew $\mathcal{V}$-category structures extending $(\mathcal{A},\underline{\mathcal{A}}(-,-))$ such that the isomorphisms \emph{(\ref{homiso})} are natural in each variable.
\end{enumerate}
\begin{proof}
By a standard argument, there are unique functors $\ast \colon \mathcal{V} \times \mathcal{A} \lra \mathcal{A}$ and \linebreak $\underline{\mathcal{A}}(-,-) \colon \mathcal{A}^\mathrm{op}\times \mathcal{A} \lra \mathcal{V}$ making the isomorphisms (\ref{actiso}) and (\ref{homiso}) natural in each variable.

(a) By the representations (\ref{actiso}), the associativity constraint for a skew $\mathcal{V}$-proactegory structure extending $(\mathcal{A},T)$ is a natural transformation
\begin{equation*}
\alpha \colon \int^B \mathcal{A}(X\ast B,C) \times \mathcal{A}(Y \ast A,B) \lra \int^Z \mathcal{A}(Z \ast A,C) \times \mathcal{V}(X\otimes Y,Z)
\end{equation*}
which, by the coend form of the Yoneda lemma, corresponds to a natural transformation
\begin{equation*}
\mathcal{A}(X \ast (Y\ast A),C) \lra \mathcal{A}((X\otimes Y)\ast A,C)
\end{equation*}
which, by the Yoneda lemma, corresponds to a natural transformation 
\begin{equation*}
a \colon (X\otimes Y) \ast A \lra X \ast (Y \ast A).
\end{equation*}
Similarly, the unit constraint for such a skew $\mathcal{V}$-proactegory structure is a natural transformation
\begin{equation*}
\lambda \colon \mathcal{A}(A,B) \lra \int^X \mathcal{A}(X \ast A,B) \times \mathcal{V}(I,X) \cong \mathcal{A}(I \ast A,B)
\end{equation*}
which, by the Yoneda lemma, corresponds to a natural transformation $l \colon I \ast A \lra A$. Moreover, the skew $\mathcal{V}$-proactegory axioms are easily shown to be equivalent under these correspondences to the skew $\mathcal{V}$-actegory axioms.

(b) By the representations (\ref{homiso}), the associativity constraint for a skew $\mathcal{V}$-proactegory structure extending $(\mathcal{A},T)$ is a natural transformation
\begin{equation*}
\alpha \colon \int^B\mathcal{V}(X,\underline{\mathcal{A}}(B,C)) \times \mathcal{V}(Y,\underline{\mathcal{A}}(A,B)) \lra \int^Z \mathcal{V}(Z,\underline{\mathcal{A}}(A,C)) \times \mathcal{V}(X \otimes Y,Z),
\end{equation*}
which, by the univeral property of coends and the coend form of the Yoneda lemma, corresponds to a natural transformation
\begin{equation*}
(\mathcal{V}\times\mathcal{V})((X,Y),(\underline{\mathcal{A}}(B,C),\underline{\mathcal{A}}(A,B)) = \mathcal{V}(X,\underline{\mathcal{A}}(B,C)) \times \mathcal{V}(Y,\underline{\mathcal{A}}(A,B)) \lra \mathcal{V}(X \otimes Y, \mathcal{A}(A,C)),
\end{equation*}
which, by the Yoneda lemma, corresponds to a natural transformation 
\begin{equation*}
M \colon \underline{\mathcal{A}}(B,C) \otimes \underline{\mathcal{A}}(A,B) \lra \underline{\mathcal{A}}(A,C).
\end{equation*}
Similarly, the unit constraint for such a skew $\mathcal{V}$-proactegory structure is a natural transformation
\begin{equation*}
\lambda \colon \mathcal{A}(A,B) \lra \int^X \mathcal{V}(X,\underline{\mathcal{A}}(A,B)) \times \mathcal{V}(I,X) \cong \mathcal{V}(I,\underline{\mathcal{A}}(A,B))
\end{equation*}
which, by the Yoneda lemma, corresponds to a natural transformation $j \colon I \lra \underline{\mathcal{A}}(A,A)$. Moreover, the skew $\mathcal{V}$-proactegory axioms are easily shown to be equivalent under these correspondences to the skew $\mathcal{V}$-category axioms.
\end{proof}
\end{proposition}

\begin{corollary} \label{yonedacor} Let $\mathcal{V}$ be a skew-monoidal category, $\mathcal{A}$ a category, and $\ast \colon \mathcal{V} \times \mathcal{A} \lra \mathcal{A}$ and $\underline{\mathcal{A}}(-,-) \colon \mathcal{A}^\mathrm{op} \times \mathcal{A} \lra \mathcal{V}$ functors. Given a natural isomorphism $\mathcal{A}(X \ast A,B) \cong \mathcal{V}(X,\underline{\mathcal{A}}(A,B))$, there is a bijection between skew $\mathcal{V}$-actegory structures extending $(\mathcal{A},\ast)$ and skew $\mathcal{V}$-category structures extending $(\mathcal{A},\underline{\mathcal{A}}(-,-))$.

Moreover, if $\mathcal{V}$ is monoidal closed, then there is a bijection between $\mathcal{V}$-actegory structures extending $(\mathcal{A},\ast)$ and tensored $\mathcal{V}$-category structures extending $(\mathcal{A},\underline{\mathcal{A}}(-,-))$.
\begin{proof}
Both kinds of structure in the first statement are equivalent, by Proposition \ref{yonedaprop}, to skew $\mathcal{V}$-proactegory structures $(\mathcal{A},T)$ for which the functors $T(X,A;-) \colon \mathcal{A} \lra \Set$ and $T(-,A;B) \colon \mathcal{V}^\mathrm{op} \allowbreak \lra \Set$ are representable for all $X \in \mathcal{V}$ and $A,B \in \mathcal{A}$. 

By the proof of Proposition \ref{yonedaprop}, the invertibility of the unit constraint $\lambda$ for such a skew $\mathcal{V}$-proactegory is equivalent both to the invertibility of the unit constraint $l$ for the corresponding skew $\mathcal{V}$-actegory, and to the normality of the corresponding skew $\mathcal{V}$-category. Similarly, the invertibility of the associativity constraints $\alpha$ and $a$ of the corresponding skew $\mathcal{V}$-proactegory and skew $\mathcal{V}$-actegory are  equivalent. Moreover, when $\mathcal{V}$ is closed, $\alpha$ corresponds to a natural transformation
\begin{equation*}
\int^B \mathcal{V}(X,\underline{\mathcal{A}}(B,C)) \times \mathcal{A}(Y \ast A,B) \lra \int^Z \mathcal{V}(Z,\underline{\mathcal{A}}(A,C)) \times \mathcal{V}(X,[Y,Z])
\end{equation*}
and hence, by the coend form of the Yoneda lemma, corresponds to a natural transformation
\begin{equation*}
\mathcal{V}(X,\underline{\mathcal{A}}(Y\ast A,C)) \lra \mathcal{V}(X,[Y,\underline{\mathcal{A}}(A,C)])
\end{equation*}
which, by the Yoneda lemma, corresponds to a natural transformation
\begin{equation*}
k \colon \underline{\mathcal{A}}(Y \ast A,C) \lra [Y,\underline{\mathcal{A}}(A,C)].
\end{equation*}
Hence the invertibility of $\alpha$ is equivalent to the invertibility of $k$.  

When $\mathcal{V}$ is monoidal closed, we may conclude that $\mathcal{V}$-actegory structures extending $(\mathcal{A},\ast)$ and tensored $\mathcal{V}$-category structures extending $(\mathcal{A},\underline{\mathcal{A}}(-,-))$ are both equivalent to $\mathcal{V}$-proactegory structures on $\mathcal{A}$ for which the functors $T(X,A;-)$ and $T(-,A;B)$ are representable.
\end{proof}
\end{corollary}

Having obtained the desired equivalence of structures, we now proceed to show that it extends to a $2$-equivalence of $2$-categories (over $\Cat$). 

\begin{definition}
A \emph{proactive functor}  $(F,\varphi) \colon (\mathcal{A},T,\alpha,\lambda) \lra (\mathcal{B},S,\alpha,\lambda)$ between skew $\mathcal{V}$-proactegories consists of the following data:
\begin{enumerate}[(i)]
\item a functor $F \colon \mathcal{A} \lra \mathcal{B}$,
\item a natural transformation $\varphi \colon T(X,A;B) \lra S(X,FA;FB)$,
\end{enumerate}
subject to the following two axioms.
\begin{equation*}
\cd{
\int^B T(X,B;C) \times T(Y,A;B) \ar[r]^-{\alpha} \ar[d]_-{\int^B \varphi \times \varphi} & \int^Z T(Z,A;C) \times P(X,Y;Z) \ar[d]^-{\int^Z \varphi \times 1} \\
\int^B S(X,FB;FC) \times S(Y,FA;FB) \ar[r]_-{\alpha} & \int^Z S(Z,FA;FC) \times P(X,Y;Z)
}
\end{equation*}
\begin{equation*}
\cd{
\mathcal{A}(A,B) \ar[r]^-{\lambda} \ar[d]_-F & \int^X T(X,A;B) \times JX \ar[d]^-{\int^X \varphi \times 1} \\
\mathcal{B}(FA,FB) \ar[r]_-{\lambda} & \int^X S(X,FA;FB) \times JX
}
\end{equation*}
\end{definition}
\begin{definition}
A \emph{proactive natural transformation} $\theta \colon (F,\varphi) \lra (G,\psi)$ between proactive functors consists of the following data:
\begin{enumerate}[(i)]
\item a natural transformation $\theta \colon F \lra G$,
\end{enumerate}
subject to the following axiom.
\begin{equation*}
\cd[@C=3.5em]{
T(X,A;B) \ar[r]^-{\varphi} \ar[d]_-{\psi} & S(X,FA;FB) \ar[d]^-{S(1,1;\theta_B)} \\
S(X,GA;GB) \ar[r]_-{S(1,\theta_A;1)} & S(X,FA;GB)
}
\end{equation*}
\end{definition}

For each skew-promonoidal category $\mathcal{V}$, skew $\mathcal{V}$-proactegories, proactive functors, and proactive natural transformations form a $2$-category $\vskproact$, with the obvious compositions. 

For each skew-monoidal category $\mathcal{V}$, let $\vskact$ denote the $2$-category of skew $\mathcal{V}$-actegories defined in \cite{MR3370862}, and let $\vact$ denote its full sub-$2$-category on the $\mathcal{V}$-actegories.

\begin{proposition} \label{yonedaprop2} Let $\mathcal{V}$ be a skew-monoidal category. Then the following are true.
\begin{enumerate}[leftmargin=*, font=\normalfont, label=(\alph*)]
\item The $2$-category $\vskact$ is $2$-equivalent over $\Cat$ to the full sub-$2$-category of \linebreak$\vskproact$ on those skew $\mathcal{V}$-proactegories $(\mathcal{A},T)$ for which the functor \linebreak $T(X,A;-) \colon \allowbreak \mathcal{A} \lra \Set$ is representable for all $X \in \mathcal{V}, A \in \mathcal{A}$.
\item The $2$-category $\vskcat$ is $2$-equivalent over $\Cat$ to the full sub-$2$-category of \linebreak$\vskproact$ on those skew $\mathcal{V}$-proactegories $(\mathcal{A},T)$ for which the functor \linebreak $T(-,A;B) \colon \allowbreak\mathcal{V}^\mathrm{op} \lra \Set$ is representable for all $A,B \in \mathcal{A}$.
\end{enumerate}
\begin{proof}
(a) Proactive functors $(F,\varphi) \colon (\mathcal{A},T) \lra (\mathcal{B},S)$ between such skew $\mathcal{V}$-proactegories correspond precisely to morphisms of skew $\mathcal{V}$-actegories under the correspondence of Proposition \ref{yonedaprop}. For the natural transformation $\varphi$ corresponds to a natural transformation
\begin{equation*}
\mathcal{A}(X\ast A,B) \lra \mathcal{B}(X\ast FA,FB)
\end{equation*} 
which corresponds, by the Yoneda lemma, to a natural transformation
\begin{equation*}
X \ast FA \lra F(X\ast A).
\end{equation*}
Under this correspondence the proactive functor axioms are equivalent to the skew $\mathcal{V}$-actegory morphism axioms.

Furthermore, a natural transformation $\theta \colon (F,\varphi) \lra (G,\psi)$ between proactive functors is proactive if and only if the diagrams
\begin{equation*}
\cd[@C=3.2em]{
\mathcal{A}(X \ast A,B) \ar[r]^-{\varphi} \ar[d]_-{\psi} & \mathcal{B}(X \ast FA,FB) \ar[d]^-{\mathcal{B}(1,\theta_B)} \\
\mathcal{B}(X\ast GA,GB) \ar[r]_-{\mathcal{B}(1\ast \theta_A,1)} & \mathcal{B}(X\ast FA,GB)
}
\end{equation*}
commute, which is equivalent by the Yoneda lemma to $\theta$ being a transformation of skew $\mathcal{V}$-actegories between morphisms of skew $\mathcal{V}$-actegories.

(b) Similarly, proactive functors $(F,\varphi) \colon (\mathcal{A},T) \lra (\mathcal{B},S)$ between such skew $\mathcal{V}$-proactegories correspond precisely to $\mathcal{V}$-functors under the correspondence of Proposition \ref{yonedaprop}. For the natural transformation $\varphi$ corresponds to  a natural transformation
\begin{equation*}
\mathcal{V}(X,\underline{\mathcal{A}}(A,B)) \lra \mathcal{V}(X,\underline{\mathcal{B}}(FA,FB))
\end{equation*}
which corresponds by the Yoneda lemma to a natural transformation 
\begin{equation*} \underline{\mathcal{A}}(A,B) \lra \underline{\mathcal{B}}(FA,FB). \end{equation*}
Under this correspondence the proactive functor axioms are equivalent to the $\mathcal{V}$-functor axioms. 

Furthermore, a natural transformation $\theta \colon (F,\varphi) \lra (G,\psi)$ between proactive functors is proactive if and only if the diagrams
\begin{equation*}
\cd[@C=4.2em]{
\mathcal{V}(X,\underline{\mathcal{A}}(A,B)) \ar[r]^-{\varphi} \ar[d]_-{\psi} & \mathcal{V}(X,\underline{\mathcal{B}}(FA,FB)) \ar[d]^-{\mathcal{V}(1,\underline{\mathcal{B}}(1,\theta_B))} \\
\mathcal{V}(X,\underline{\mathcal{B}}(GA,GB)) \ar[r]_-{\mathcal{V}(1,\underline{\mathcal{B}}(\theta_A,1))} & \mathcal{V}(X,\underline{\mathcal{B}}(FA,GB))
}
\end{equation*}
commute, which is equivalent by the Yoneda lemma to $\mathcal{V}$-naturality of $\theta$ between the corresponding $\mathcal{V}$-functors.
\end{proof}
\end{proposition}

\begin{corollary}
Let $\mathcal{V}$ be a skew-monoidal category. Then the full sub-$2$-category of $\vskact$ on those skew $\mathcal{V}$-actegories $(\mathcal{A},\ast)$ for which $-\ast A \colon \mathcal{V} \lra \mathcal{A}$ has a right adjoint for all $A \in \mathcal{A}$ is $2$-equivalent over $\Cat$ to the full sub-$2$-category of $\vskcat$ on those skew $\mathcal{V}$-categories $(\mathcal{A},\underline{\mathcal{A}}(-,-))$ for which $\underline{\mathcal{A}}(A,-) \colon \mathcal{A} \lra \mathcal{V}$ has a left adjoint for all $A \in \mathcal{A}$.

Moreover, if $\mathcal{V}$ is monoidal closed, then the full sub-$2$-category of $\vact$ on those $\mathcal{V}$-actegories $(\mathcal{A},\ast)$ for which $-\ast A \colon \mathcal{V} \lra \mathcal{A}$ has a right adjoint for all $A \in \mathcal{A}$ is $2$-equivalent over $\Cat$  to the full sub-$2$-category of $\vcat$ on the tensored $\mathcal{V}$-categories.
\begin{proof}
These results follow from Proposition \ref{yonedaprop2} as Corollary \ref{yonedacor} follows from Proposition \ref{yonedaprop}.
\end{proof}
\end{corollary}

\bibliographystyle{alpha}

\end{document}